\documentclass[12pt, leqno]{amsart}
\usepackage{array}
\usepackage{longtable}
\usepackage{url}
\usepackage{enumerate}
\usepackage{amsmath}
\usepackage{amssymb}
\usepackage{amsthm}
\usepackage{amsfonts}
\usepackage{tikz,pgf}
\usepackage{graphicx}
\usepackage{tikz-cd}
\usepackage{mathtools}
\usepackage{tikz}
\usetikzlibrary{arrows}
  \usepackage{ifthen}
 \newcommand{\mymarginpar}[1]{%
    \marginpar{\ifthenelse{\isodd{\arabic{page}}}{\flushleft 
#1}{\flushright #1}}}

\numberwithin{equation}{section}

 \newcommand{\eps}{\varepsilon}           
 
 \newcommand{\IC}{\mathbb{C}}

 \newcommand{\IN}{\mathbb{N}}                  
                   
 \newcommand{\IQ}{\mathbb{Q}}                  
 \newcommand{\IR}{\mathbb{R}}                  
                   
 \newcommand{\IT}{\mathbb{T}}                  
 \newcommand{\IZ}{\mathbb{Z}}

\newcommand{\CK}{\mathcal{K}}

 \theoremstyle{plain} 
 \newtheorem{Theorem}{Theorem}[section]
 \newtheorem{Lemma}[Theorem]{Lemma}
 \newtheorem{Proposition}[Theorem]{Proposition}
 
 \newtheorem{Corollary}[Theorem]{Corollary}

 \theoremstyle{definition} 
 \newtheorem{Definition}[Theorem]{Definition}
 \newtheorem{Remark}[Theorem]{Remark}

 \newtheorem{Notation}[Theorem]{Notation}

\allowdisplaybreaks[1]

\begin{document}

\title[$C^*$-algebra extensions and rational numbers]{$C^*$-algebra extensions associated to continued fraction expansions of rational numbers}


\author[Jack Spielberg]{Jack Spielberg}

\address{School of Mathematical and Statistical Sciences\\Arizona State University\\USA}
\email{jack.spielberg@asu.edu}

\begin{abstract}

We solve the isomorphism problem for essential unital $C^*$-algebra extensions of the form $0 \to \CK \oplus \CK \to E \xrightarrow{\pi} M_n \otimes C(\IT) \to 0$. We then relate these to analogs of the Effros-Shen AF algebras for rational numbers. This involves $C^*$-algebras constructed from categories of paths built from certain nonsimple continued fraction expansions.

\end{abstract}

\maketitle

\section{Introduction}
\label{sec intro}

For an irrational number $\theta \in (0,1)$ there corresponds a simple unital AF algebra $A_\theta$ with Elliott invariant $(K_0(A_\theta),K_0(A_\theta)_+,[1]_0) = (\IZ + \theta \IZ,(\IZ + \theta \IZ) \cap [0,\infty),1)$. In \cite{es}, Effros and Shen construct the algebra $A_\theta$ explicitly from the simple continued fraction expansion of $\theta$. They use the rational approximants of the expansion to define an inductive sequence of finite dimensional $C^*$-algebras whose limit is $A_\theta$.
In \cite{ms} an alternative construction of $A_\theta$ is given, based on a certain nonsimple continued fraction expansion of $\theta$. In this version $A_\theta$ is given as an inductive limit of infinite dimensional type $I$ $C^*$-algebras. This construction has some surprising features; for example, $A_\theta$ is given as the $C^*$-algebra of a non-AF groupoid.

A basic feature of continued fractions is that while an irrational number has a unique (infinite) simple continued fraction expansion, each rational number has two distinct finite simple continued fraction expansions. In the setting of \cite{es} this means that there are two finite dimensional $C^*$-algebras that can play the role of the ``Effros-Shen algebra'' for a rational number.
It is a remarkable feature of the construction in \cite{ms} that the continued fraction expansions used have only one version for rational numbers as well as for irrational numbers. This provides a unique candidate for an analog of the Effros-Shen algebras at a rational number. This candidate is a certain $C^*$-algebra extension of a matrix algebra over $C(\IT)$ by the direct sum of two copies of the compact operators.

In this paper we classify these extension algebras up to isomorphism, thus identifying the class of ``pseudo Effros-Shen algebras'' for rational numbers that are mentioned above. In future work we will show that when combined with the usual Effros-Shen algebras for irrational numbers, we obtain an upper semicontinuous field over the interval $[0,1)$. The extension classification we treat is nonstandard. The usual classification of $C^*$-algebra extensions is by means of the \textit{Busby invariant}, up to a suitable notion of equivalence. Our problem, rather, is to classify the extension algebra itself up to $*$-isomorphism. In sections \ref{sec extensions} and \ref{sec the invariant} we solve this extension problem in general. The invariant involves the index of the extension, and the defect numbers of the two ancillary extensions of the matrix algebra by the compact operators. We also specialize to the case corresponding to our analogs of Effros-Shen algebras at rational numbers, completing the classification of such algebras. In section \ref{sec continued fractions of rationals} we prove that there is a one-to-one correspondence between the data defining our categories of paths and a class of the extension algebras studied in the previous sections. For this we rely on the work in \cite{ms}. Therefore we begin the section with a very brief description of the construction of a $C^*$-algebra from a rational number along the lines of that paper. This construction begins with a category of paths built from a certain nonsimple continued fraction expansion of a number in $[0,1)$, and then from the groupoid and $C^*$-algebra associated to the category as in \cite{spi}. The details of this construction are needed in only a few places, so we have tried to give just enough background to make those parts intelligible. We have avoided giving an introduction to categories of paths, and the groupoids and $C^*$-algebras associated to them. The correspondence between a rational number and the invariant of an extension algebra appears to be of independent interest. It involves a variation on the usual Euclidean algorithm for finding the greatest common divisor of two integers. In the final section \ref{sec continued fractions} we discuss the nonsimple continued fraction expansions that correspond to the construction, here and in \cite{ms}, of categories of paths from numbers in $[0,1)$. Then we can state our result on the bijection between rational numbers and the isomorphism classes of extension $C^*$-algebras studied in the earlier sections.

We let $\IN$ denote the nonegative integers $\{0,1,2,\ldots\}$.

\section{Classifying certain $C^*$-algebra extensions}
\label{sec extensions}

Let $H_+$ and $H_-$ be two copies of a separable infinite dimensional Hilbert space. We will classify essential unital extensions
\begin{equation} \label{eqn the extension}
0 \to \CK(H_+) \oplus \CK(H_-) \to E \xrightarrow{\pi} M_n \otimes C(\IT) \to 0
\end{equation}
up to isomorphism of the algebra $E$. We will usually write $\CK \oplus \CK$ for the ideal. Since $E$ is a subalgebra of $B(H_+) \oplus B(H_-)$, each element $x$ of $E$ has the form $x = x_+ \oplus x_-$, with $x_{\pm} \in B(H_\pm)$ (and we may write $+/-$ as superscripts occasionally). The usual study of extensions classifies the \textit{Busby invariant} $\tau : M_n \otimes C(\IT) \hookrightarrow Q(H_+) \oplus Q(H_-)$, where $Q(H) = B(H)/\CK(H)$ is the Calkin algebra, up to one of several equivalence relations (\cite[section 15.2]{blackadar}). Our problem is different, but we will make use of the Busby invariant when we find it convenient.

The $C^*$-algebra $M_n \otimes C(\IT)$ is generated in the usual way by matrix units $e_{i1} \otimes 1$, $2 \le i \le n$, and a partial unitary $e_{11} \otimes z$. Thus any $*$-homomorphism from $M_n \otimes C(\IT)$ is determined by its values on these elements. There is an automorphism $\eta$ of $M_n \otimes C(\IT)$ defined by $\eta(e_{i1} \otimes 1) = e_{i1} \otimes 1$ and $\eta(e_{11} \otimes z) = e_{11} \otimes z^{-1}$. Then $\tau$ and $\tau \circ \eta$ define the same algebra $E$. Since $H_+ = H_-$, there is an automorphism $\eps$ of $B(H_+) \oplus B(H_-)$ given by $\eps(x_+ \oplus x_-) = x_- \oplus x_+$, and $\eps$ descends to an automorphism of $Q(H_+) \oplus Q(H_-)$. The Busby invariant $\eps \circ \tau$ determines the algebra $\pi^{-1} ( \eps \circ \tau(M_n \otimes C(\IT)) ) = \eps(E)$. Of course, $\eps(E) \cong E$.

The first isomorphism invariant of $E$ comes from the index map in $K$-theory, or equivalently, from the Fredholm index.

\begin{Definition} \label{def index}

Given the extension \eqref{eqn the extension}, the \emph{index} of $\tau$ is $\text{ind}(\tau) := \partial([e_{11} \otimes z]_1) \in K_0(\CK \oplus \CK) \cong \IZ^2$.

\end{Definition}

\begin{Definition} \label{def matrix units}

Let $A$ be a $C^*$-algebra. A \textit{family of $n \times n$ matrix units} in $A$ is a set $\{ e_{ij} : 1 \le i,j \le n \}$ such that $e_{ij} = e_{ji}^*$ and $e_{ij} e_{k \ell} = \delta_{jk} e_{i\ell}$ for all $i,j,k,\ell$.

\end{Definition}

\begin{Remark} \label{rem lifting matrix units}

It is a fairly routine exercise to use \cite[Theorem 2.5]{calkin} to prove the following: if $\{ e_{ij} \}$ is a family of $n \times n$ matrix units in $Q(H)$, then there is a family, $\{ E_{ij} \} \subseteq B(H)$ of $n \times n$ matrix units such that $\pi(E_{ij}) = e_{ij}$ for all $i,j$.

\end{Remark}

If $\{E_{ij}\} \subseteq E$ are matrix units lifting $\{e_{ij} \otimes 1\}$, and if $S \in E_{11} E E_{11}$ lifts $e_{11} \otimes z$, then $S_\pm$ is a Fredholm operator in $B(E_{11}^\pm H_\pm)$ with Fredholm index $a_\pm$, where ind$(\tau) = (a_+,a_-)$.

It is clear that ind$(\tau \circ \eta) = - \text{ind}(\tau)$ and that ind$(\eps \circ \tau) = \eps(\text{ind}(\tau))$, where we write $\eps$ also for the interchange of coordinates on $\IZ^2$. Thus the isomorphism class of $E$ distinguishes only the orbit of ind$(\tau)$ under the action of $(-1)\cdot$ and $\eps$.

\begin{Definition} \label{def m}

For a family of $n \times n$ matrix units $\{E_{ij}\} \subseteq B(H)$, the \emph{defect} of $\{E_{ij}\}$ is the rank of the projection $1 - \sum_{i=1}^n E_{ii} \in \IN \cup \{\infty\}$. For a family of $n \times n$ matrix units $\{E_{ij}\} \subseteq B(H_+) \oplus B(H_-)$,  we write $m(\{E_{ij}\}) = (k_+,k_-)$, where $k_\pm$ is the defect of $\{E^\pm_{ij}\}$ in $B(H_\pm)$.

\end{Definition}

Suppose that $\{ E_{ij} \}$ are lifts of $\{ e_{ij} \otimes 1 \}$ in \eqref{eqn the extension}. Then $\pi(1 - \sum_i E_{ii}) = 1 - \sum_i e_{ii} \otimes 1 = 0$, so $1 - \sum_i E_{ii} \in \CK \oplus \CK$. Therefore the defect of $1^\pm - \sum_i E_{ii}^\pm$ is finite, so $m(\{ E_{ij}) \in \IN^2$.

\begin{Remark} \label{rem unitary equivalence of matrix units}

Let $H$ be an infinite dimensional Hilbert space. If $\{ E_{ij} \}$ and $\{ F_{ij} \}$ are two families of $n \times n$ matrix units in $B(H)$ with the same finite defect $k$, then there is $W \in U(H)$ such that $W E_{ij} W^* = F_{ij}$ for all $i,j$. To see this, note that $E_{ii} H$ and $F_{ii} H$ are all infinite dimensional. Let $W_0 \in U(E_{11}H,F_{11}H)$ and $W_1 \in U( (1 - \sum_i E_{ii}) H, (1 - \sum_i F_{ii}) H)$ be arbitrary. Then we may let $W = \sum_i F_{i1} W_0 E_{1i} + W_1$.

\end{Remark}

\begin{Lemma} \label{lem lift the partial unitary}

Let \eqref{eqn the extension} be an extension, and let $\{E_{ij} : 1 \le i,j \le n \} \subseteq E$ be matrix units with $\pi(E_{ij}) = e_{ij} \otimes 1$. Let \emph{ind}$(\tau) = (a_+,a_-)$. There is a partial isometry $S = S_+ \oplus S_- \in B(E_{11}^+ H_+) \oplus B(E_{11}^- H_-)$ such that $\pi(S) = e_{11} \otimes z$, and such that $S_\pm$ is an isometry, respectively a coisometry, if $a_\pm \le 0$, respectively if $a_\pm \ge 0$. We note that $E$ is generated as a $C^*$-algebra by $\{ E_{ij} \}$, $S$, and $\CK \oplus \CK$.

\end{Lemma}

\begin{proof}
Let $u = e_{11} \otimes z = u_+ \oplus u_- \in Q(E_{11}^+ H_+) \oplus Q(E_{11}^- H_-)$, with $u_\pm$ unitary. The lemma now follows from \cite[Theorem 2.7]{calkin}.
\end{proof}

\begin{Lemma} \label{lem shift defect by index}

Let \eqref{eqn the extension} be an extension with Busby invariant $\tau$. Let $\{E_{ij} : 1 \le i,j \le n \} \subseteq E$ be matrix units with $\pi(E_{ij}) = e_{ij} \otimes 1$. Let \emph{ind}$(\tau) = (a_+,a_-)$. Suppose that $m(\{E_{ij}\}) - (a_+,a_-) \ge 0$, where this refers to the product order on $\IZ^2$. Then there are matrix units $\{F_{ij}\} \subseteq E$ such that $\pi(F_{i1}) = e_{i1} \otimes z^{-1}$ for $2 \le i \le n$, $\pi(F_{11}) = e_{11} \otimes 1 = \pi(E_{11})$, and $m(\{F_{ij}\}) = m(\{E_{ij}\}) + (a_+,a_-)$.

\end{Lemma}

\begin{proof}
Let $m(\{E_{ij}\}) = (k_+,k_-)$, and let $S = S_+ \oplus S_- \in B(E_{11}^+ H_+) \oplus B(E_{11}^- H_-)$ be as in Lemma \ref{lem lift the partial unitary}. Suppose for definiteness that $a_+ < 0$ and $a_- > 0$. Let $\{h_{\alpha \beta}\}_{\alpha,\beta \ge 1}$ be the usual matrix units in $B(\ell^2)$ and let $R$ be the right shift on $\ell^2$ (so that $R h_{\alpha \beta} = h_{\alpha + 1, \beta + 1} R$ and $h_{11}R = 0$). By the Wold decomposition we may write
\begin{align*}
H_\pm &= \bigl[ \IC^n \otimes \bigl( (\ell^2 \otimes \IC^{|a_\pm|}) \oplus L_\pm \bigr) \bigr] \oplus \IC^{k_\pm}, \\
\noalign{\noindent for Hilbert spaces $L_\pm$, and}
 S_\pm &= \bigl[ e_{11} \otimes \bigl( (R^{\pm 1} \otimes 1) \oplus V_\pm \bigr) \bigr] \oplus 0,
\end{align*}
where $E_{ij}^\pm = \bigl[ e_{ij} \otimes  \bigl( (1 \otimes 1) \oplus 1 \bigr) \bigr] \oplus 0$, $V_\pm \in U(L_{\pm})$, and we let $R^{-1}$ denote $R^*$. (Note that since we are considering the case where $a_+ < 0$, $S_+$ is an isometry, so we use $R^{+1}$ to model $S_+$, and similarly we use $R^{-1} = R^*$ to model the coisometry $S_-$.) Since $k_- - a_-  \ge 0$ there is a projection $q'$ of rank $a_-$ in $B(\IC^{k_-}) \subseteq B(H_-)$. Let $q = \bigl[ e_{11} \otimes \bigl( h_{11} \otimes 1) \oplus 0 \bigr) \bigr] \oplus 0 \in B(H_-)$; $q$ is also a projection of rank $a_-$. Let $w$ be  a partial isometry with $w^* w = q$ and $w w^* = q'$. Let $p = \bigl[ e_{11} \otimes \bigl( h_{11} \otimes 1) \oplus 0 \bigr) \bigr] \oplus 0 \in B(H_+)$; $p$ is a projection of rank $|a_+|$.

Now we set
\begin{align*}
F_{21} &= E_{21} (S^* + w^*), \\
\noalign{\noindent and for $2 < i \le n$,}
F_{i1} &= E_{i2} F_{21}.
\end{align*}
Then $\pi(F_{21}) = (e_{21} \otimes 1)(e_{11} \otimes z)^* = e_{21} \otimes z^{-1}$, and for $i > 2$ we have $\pi(F_{i 1}) = (e_{i 2} \otimes 1)(e_{21} \otimes z^{-1}) = e_{i 1} \otimes z^{-1}$. We claim that
\begin{align}
F_{i 1} F_{i 1}^* &= E_{i i}, \text{ for $2 \le i \le n$,} \label{lem shift defect by index one} \\
F_{i 1}^* F_{i 1} &= (E_{11}^+ - p) \oplus (E_{11}^- + q'), \text{ for $2 \le i \le n$.} \label{lem shift defect by index two}
\end{align}
%
%
To see this, first note that
%
\begin{align}
w^* S &= (w^* q' )(E_{11} S) = w^* (q' E_{11}) S = 0, \label{lem shift defect by index three} \\
Sw^*& = S q w^* = 0_{H_+} \oplus \Bigl( \bigl[ e_{11} \otimes \bigl( (R^* h_{11} \otimes 1) \oplus 0 \bigr) \oplus 0 \bigr] \Bigr) w^* =  0. \label{lem shift defect by index four}
\end{align}
%

Now we prove \eqref{lem shift defect by index one}. 
\begin{align*}
F_{21} F_{21}^*  &= E_{21} (S^* + w^*)(S + w) E_{12} \\
&= E_{21} (S^*S + w^* w) E_{12}, \text{ by \eqref{lem shift defect by index three},} \\
&= E_{21} (S^*S + q) E_{12} \\
&= E_{21} \Bigl(E_{11}^+ \oplus \bigl(E_{11}^- - \bigl( \bigl[ e_{11} \otimes \bigl( (h_{11} \otimes 1) \oplus 0 \bigr) \bigr] \oplus 0 \bigr) + q \bigr) \Bigr) E_{12} \\
&= E_{21} (E_{11}^+ \oplus E_{11}^-) E_{12} \\
&= E_{22}.
\end{align*}
For $i > 2$ we have $F_{i1} F_{i1}^* = E_{i2} F_{21} F_{21}^* E_{2i} = E_{i2} E_{22} E_{2i} = E_{ii}$. To prove \eqref{lem shift defect by index two}, we have
\begin{align*}
F_{21}^* F_{21} &= (S + w) E_{12} E_{21} (S^* + w^*) \\
&= SS^* + w w^*, \\
\noalign{\noindent by \eqref{lem shift defect by index four}, and the fact that $wE_{11} = wqE_{11} = wq = w$,}
&= SS^* + q' \\
&= (E_{11}^+ - p) \oplus (E_{11}^- + q').
\end{align*}
For $i > 2$ we have $F_{i1}^* F_{i1} = F_{21}^* E_{i2}^* E_{i2} F_{21} = F_{21}^* E_{22} F_{21} = F_{21}^* F_{21}$, by \eqref{lem shift defect by index one}.

It follows that the $\{F_{i1}^* F_{i1}\}$ are all equal, and that $\{F_{21}^* F_{21} \} \cup \{F_{i1} F_{i1}^* : 2 < i \}$ are pairwise orthogonal projections. Therefore the $\{F_{i1}\}$ generate a family of $n\times n$ matrix units in $E$, and $m(\{F_{ij}\}) = (k_+ - |a_+|, k_- + a_-) = (k_+ + a_+, k_- + a_-)  = m(\{E_{ij}\}) + (a_+,a_-)$.
\end{proof}

\begin{Remark} \label{rem shift defect by index}
\begin{enumerate}

\item \label{rem shift defect by index one} Note that with $S$ as in the proof of Lemma \ref{lem shift defect by index} we have that $E$ is generated as a $C^*$-algebra by $\{ F_{ij}  \}$, $S$, and $\CK \oplus \CK$.

\item \label{rem shift defect by index two} If $m(\{E_{ij}\}) + (a_+,a_-) \ge 0$ then we may alter the proof by choosing a projection $p' \in B(\IC^{k_+})$ of rank $|a_+|$, and a partial isometry $v$ with $v^* v = p$ and $vv^* = p'$, and setting $F_{21} = E_{21} S + v^*$ (and then $\pi(F_{i1}) = e_{11} \otimes z$, $i \ge 2$).

\item \label{rem shift defect by index three} It is straightforward to modify the proof for any configuration of the signs of $a_+$ and $a_-$.

\end{enumerate}
\end{Remark}

\begin{Lemma} \label{lem adjust by n}

Let an extension \eqref{eqn the extension} be given, and let $\{E_{ij}\}$ be a family of matrix units lifting $\{e_{ij} \otimes 1\}$. Let $\mu_\pm \in \IN$. Then there is a family of matrix units $\{E_{ij}'\}$ also lifting $\{e_{ij} \otimes 1\}$ such that $m(\{E_{ij}'\}) = m(\{E_{ij}\}) + n(\mu_+,\mu_-)$. Moreover $C^* \bigl( \{E_{ij}'\}, \CK \oplus \CK \bigr) = C^* \bigl( \{E_{ij}\}, \CK \oplus \CK \bigr)$, 

\end{Lemma}

\begin{proof}
Choose projections $c_\pm \le E_{11}^\pm$ with rank$(c_\pm) = n\mu_\pm$. Let $E_{11}' = E_{11} - c_+ \oplus c_-$, and set $E_{i1}' = E_{i1} E_{11}'$ for $i > 1$.
\end{proof}

\begin{Definition} \label{def invariant of extension}

Let an extension \eqref{eqn the extension} be given, and let $a = (a_+,a_-) = \text{ind}(\tau)$. Let $D = \IZ^2 / (\IZ a + n\IZ^2)$. For a family of $n \times n$ matrix units $\{F_{ij}\} \subseteq E$ we write $\overline{m}(\{F_{ij}\}) := m(\{F_{ij}\}) + \IZ a + n\IZ^2 \in D$.

\end{Definition}

The next proposition shows that $\overline{m}$ is an invariant of the algebra $E$.

\begin{Proposition} \label{prop invariant of extension}

Let an extension \eqref{eqn the extension} be given, and let $a = (a_+,a_-) = \text{\emph{ind}}(\tau)$. Let $\{E_{ij}\}$ be matrix units lifting $\{e_{ij} \otimes 1\}$, and let $\{F_{ij}\}$ be a family of $n \times n$ matrix units in $E$ with $1 - \sum_{i=1}^n F_{ii} \in \CK \oplus \CK$. Then $\overline{m}(\{E_{ij}\}) = \overline{m}(\{F_{ij}\})$.

\end{Proposition}

\begin{proof}
We will make a series of modifications to $\{F_{ij}\}$ that do not change $\overline{m}(\{F_{ij}\})$. For each $i$, $\pi(F_{ii})$ is a minimal projection in $M_n \otimes C(\IT) \cong C(\IT,M_n)$, that is, it is a rank one projection at each point of $\IT$. The space of rank one projections in $M_n$ is the complex projective space $\IC P^n$, which is well known to be connected and simply connected. Thus any two continuous maps from $\IT$ to $\IC P^n$ are homotopic. Therefore $\pi(F_{11})$ and $\pi(E_{11})$ are homotopic in $C(\IT,\IC P^n)$. It follows that there is a continuous path $t \in [0,1] \mapsto w_t \in U(M_n \otimes C(\IT))$ such that $w_0 = 1$ and $w_1(e_{11} \otimes 1)w_1^* = \pi(F_{11})$ (e.g. apply the polar part to \cite[4.3.3]{blackadar}). Recall that a unitary element $u \in Q(H)$ lifts to a unitary element in $B(H)$ if and only if $u$ has Fredholm index zero, equivalently $u$ is homotopic to 1 in $UQ(H)$. Thus there is a unitary $W_1 \in B(H_+) \oplus B(H_-)$ such that $\pi(W_1) = w_1$. Replacing $F_{ij}$ by $W_1^* F_{ij} W_1$, we may assume that $\pi(F_{11}) = e_{11} \otimes 1$.

Now consider
\begin{equation*}
\{\pi(F_{ij}) : i,j \ge 2\} \subseteq (\sum_{i=2}^n e_{ii} \otimes 1)(M_n \otimes C(\IT))(\sum_{i=2}^n e_{ii} \otimes 1) \cong M_{n-1} \otimes C(\IT).
\end{equation*}
Applying the above argument to $e_{22} \otimes 1$ and $\pi(F_{22})$, we obtain a unitary $W_2 \in (\sum_{i=2}^n F_{ii}) (B(H_+) \oplus B(H_-)) (\sum_{i=2}^n F_{ii})$ such that $\pi(W_2) (e_{22} \otimes 1) \pi(W_2)^* = \pi(F_{22})$. Replacing $F_{ij}$ by $(F_{11} + W_2)^* F_{ij} (F_{11} + W_2 )$, we may assume that $\pi(F_{ii}) = e_{ii} \otimes 1$ for $i = 1,2$. Repeating this process, we may assume that $\pi(F_{ii}) = e_{ii} \otimes 1$ for $1 \le i \le n$.

For each $i > 1$, $\pi(E_{1i}F_{i1}) \in U(e_{11} \otimes C(\IT))$. Then there is $\phi_i \in C(\IT,\IT)$ such that $\pi(E_{1i}F_{i1}) = e_{11} \otimes \phi_i$. It follows that $\pi(F_{i1}) = e_{i1} \otimes \phi_i$. Let $r_i = \omega(\phi_i)$, the winding number of $\phi_i$, for $2 \le i \le n$. Let $r_1 = 0$ and $\phi_1 = 1$. Then $\phi_i$ is homotopic to $z^{r_i}$ in $C(\IT,\IT)$. We have
\begin{equation*}
\pi(F_{ij}) = \pi(F_{i1}) \pi(F_{1j}) = (e_{i1} \otimes \phi_i)(e_{1j} \otimes \overline{\phi_j}) = e_{ij} \otimes \phi_i \overline{\phi_j}.
\end{equation*}
Let $b = \sum_{i=1}^n e_{ii} \otimes z^{r_i} \overline{\phi_i}$. Since $\omega(z^{r_i} \overline{\phi_i}) = 0$ for all $i$, $b$ is homotopic to 1 in $U(M_n \otimes C(\IT))$. Then there is a unitary $B \in B(H_+) \oplus B(H_-)$ such that $\pi(B) = b$. Note that
\begin{equation*}
\pi(B F_{ij} B^*) = (e_{ii} \otimes z^{r_i} \overline{\phi_i}) (e_{ij} \otimes \phi_i \overline{\phi_j}) (e_{jj} \otimes z^{-r_j} \phi_j) = e_{ij} \otimes z^{r_i - r_j}.
\end{equation*}
Replacing $\{F_{ij}\}$ by $\{B F_{ij} B^*\}$, we may assume that $\pi(F_{ij}) = e_{ij} \otimes z^{r_i - r_j}$. Let $s_1 = \sum_{i=2}^n r_i$, and $s_i = -r_i$ for $i > 1$. Let $w = \sum_{i=1}^n e_{ii} \otimes z^{s_i}$. Then $w \in U(M_n \otimes C(\IT)$, and ind$(w) = \partial[w]_1 = -\sum_{i=1}^n s_i = 0$. It follows that there is a unitary operator $W \in B(H_+) \oplus B(H_-)$ such that $\pi(W) = w$. Now, for $i > 1$ we have
\begin{equation*}
\pi(W F_{i1} W^*) = (e_{ii} \otimes z^{s_i})(e_{i1} \otimes z^{r_i - r_1})(e_{11} \otimes z^{-s_1}) = e_{i1} \otimes z^{-s_1}.
\end{equation*}
Replacing $\{F_{ij}\}$ by $\{W F_{ij} W^*\}$, we may assume that $\pi(F_{i1}) = e_{i1} \otimes z^{-s_1}$ for $i > 1$, and that $\pi(F_{11}) = e_{11} \otimes 1$.

We now want to apply Lemma \ref{lem shift defect by index} $|s_1|$ times. By Lemma \ref{lem adjust by n} there is a family of matrix units $\{F_{ij}'\}$ such that $\pi(F_{ij}') = \pi(F_{ij})$, $m(\{F_{ij}'\}) \in m(\{F_{ij}\}) + n \IZ^2$, and $m(\{F_{ij}'\}) \ge m(\{F_{ij}\}) + |s_1|(|a_+|,|a_-|)$. Then $\overline{m}(\{F_{ij}\}) = \overline{m}(\{F_{ij}'\})$. We may now apply Lemma \ref{lem shift defect by index} (and Remarks \ref{rem shift defect by index}) $|s_1|$ times to obtain matrix units $\{\widetilde{F}_{ij}\}$ with $\overline{m}(\{\widetilde{F}_{ij}\}) = \overline{m}(\{F_{ij}\})$ and such that $\pi(\widetilde{F}_{ij}) = e_{ij} \otimes 1$. It follows that $\{E_{ij}^\pm \}$ and $\{\widetilde{F}_{ij}^\pm \}$ define the same Busby invariants $\tau_\pm : M_n \to Q(H_\pm)$. It follows that the defects of $\{E_{ij}^\pm \}$ and $\{ \widetilde{F}_{ij}^\pm \}$ are equal modulo $n$ (\cite[Proposition 2]{thayer}). Then $\overline{m}(\{E_{ij}\}) = \overline{m}(\{\widetilde{F}_{ij}\}) = \overline{m}(\{F_{ij}\})$.
\end{proof}

\begin{Corollary} \label{cor invariant of extension}

Let $E$ and $F$ be the algebras in the center places of two versions of the extension \eqref{eqn the extension} (with corresponding Busby invariants $\tau_E$ and $\tau_F$). Suppose that $\varphi : E \to F$ is an isomorphism such that $\varphi(\CK \oplus 0) = \CK \oplus 0$ (and hence also that $\varphi(0 \oplus \CK) = 0 \oplus \CK$). Then \emph{ind}$(\tau_F) = \pm \emph{ind}(\tau_E)$. Moreover, for any families of matrix units $\{E_{ij}\} \subseteq E$ and $\{F_{ij}\} \subseteq F$ lifting $\{e_{ij} \otimes 1\}$, we have that $\overline{m}(\{E_{ij}\}) = \overline{m}(\{F_{ij}\})$.

\end{Corollary}

\begin{proof}
Without loss of generality we may as well assume that $E,F \subseteq B(H_+) \oplus B(H_-)$. The hypotheses on $\varphi$ imply that $\varphi = \varphi_+ \oplus \varphi_-$, where $\varphi_\pm |_{\CK(H_\pm)} \in Aut (\CK(H_\pm))$. It follows that there are $V_\pm \in U(H_\pm)$ such that $\varphi_\pm = Ad(V_\pm)$. Let $V = V_+ \oplus V_-$. Then $\varphi = Ad(V)$, and letting $\overset{\boldsymbol{\cdot}}{\varphi} = Ad(\pi(V))$, we have $\overset{\boldsymbol{\cdot}}{\varphi}  \circ \pi|_E = \pi|_F \circ \varphi$. Since $\overset{\boldsymbol{\cdot}}{\varphi} \in Aut(M_n \otimes C(\IT))$, $\overset{\boldsymbol{\cdot}}{\varphi}_* = \pm \text{id}$ on $\IZ \cong K_1(M_n \otimes C(\IT))$. Then $\partial_E = \partial_F \circ \overset{\boldsymbol{\cdot}}{\varphi}$, and hence ind$(\tau_F) = \pm \text{ind}(\tau_E)$. It follows that the group $D$ from Definition \ref{def invariant of extension} is the same for $E$ and $F$. Proposition \ref{prop invariant of extension} implies that $\overline{m}(\{\varphi(F_{ij})\}) = \overline{m}(\{E_{ij}\})$. Since $\varphi_\pm$ are implemented by unitaries, $m(\{F_{ij}\}) = m(\{\varphi(F_{ij})\})$. Therefore $\overline{m}(\{E_{ij}\}) = \overline{m}(\{F_{ij}\})$.
\end{proof}

\begin{Remark} \label{rem invariant of extension}

In the context of Corollary \ref{cor invariant of extension}, if the isomorphism $\varphi$ interchanges the two copies of $\CK$, then $\eps \circ \varphi : E \to \eps(F)$ is an automorphism satisfying the hypotheses of the Corollary. Since $F$ and $\eps(F)$ are isomorphic, the hypothesis on $\varphi$ in Corollary \ref{cor invariant of extension} is not a serious limitation.

\end{Remark}

\begin{Proposition} \label{prop invariant implies isomorphism}

Let $E$ and $F$ be the algebras in the center places of two versions of the extension \eqref{eqn the extension} (with corresponding Busby invariants $\tau_E$ and $\tau_F$) such that \emph{ind}$(\tau_E) = \emph{ind}(\tau_F)$.  Let $\{E_{ij}\} \subseteq E$ and $\{F_{ij}\} \subseteq F$ be families of matrix units lifting $\{e_{ij} \otimes 1\}$. Suppose that $\overline{m}(\{E_{ij}\}) = \overline{m}(\{F_{ij}\})$. Then $E$ and $F$ are isomorphic.

\end{Proposition}

\begin{proof}
Without loss of generality we may as well assume that $E,F \subseteq B(H_+) \oplus B(H_-)$. Let $a = (a_+,a_-) = \text{ind}(\tau_E) = \text{ind}(\tau_F)$. There are $\alpha,\beta_\pm \in \IZ$ such that $m(\{F_{ij}\}) = m(\{E_{ij}\}) + \alpha a +n \beta$. Choose $\nu_\pm \in \IN$ such that $m(\{E_{ij}\}) + n\nu \ge (|\alpha a_+|, |\alpha a_-|)$ and $\nu \ge \beta$. By Lemma \ref{lem adjust by n} there is a family $\{E_{ij}^{(2)}\}$ also lifting $\{e_{ij} \otimes 1\}$ with $m(\{E_{ij}^{(2)}\}) = m(\{E_{ij}\}) + n\nu$. By Lemma \ref{lem shift defect by index} and Remark \ref{rem shift defect by index} (applied $|\alpha|$ times), there is a family of matrix units $\{E_{ij}^{(3)}\} \subseteq E$ such that $m(\{E_{ij}^{(3)}\}) = m(\{E_{ij}^{(2)}\}) + \alpha a$, and such that $\pi(E_{i1}^{(3)}) = e_{i1} \otimes z^{-\alpha}$ and $\pi(E_{11}^{(3)}) = e_{11} \otimes 1$. Then
\begin{align*}
m(\{E_{ij}^{(3)}\}) &= m(\{E_{ij}\}) + \alpha a + n\nu \\
&= m(\{E_{ij}\}) + \alpha a +n\beta + n(\nu - \beta) \\
&= m(\{F_{ij}\}) + n(\nu - \beta).
\end{align*}
Since $\nu - \beta \ge 0$, Lemma \ref{lem adjust by n} implies that there is a family of matrix units $\{ F^{(2)}_{ij} \} \subseteq F$ also lifting $\{ e_{ij} \otimes 1 \}$ and with $m(\{F^{(2)}_{ij} \}) = m(\{F_{ij}\}) + n(\nu - \beta)$. Then $m(\{E_{ij}^{(3)}\}) = m(\{F^{(2)}_{ij} \})$. It follows from Remark \ref{rem unitary equivalence of matrix units} that there is a unitary $W \in B(H_+) \oplus B(H_-)$ such that $W E_{ij}^{(3)} W^* = F^{(2)}_{ij}$, $1 \le i,j \le n$. Let $S$ be as in Lemma \ref{lem lift the partial unitary} for $\{E^{(3)}_{ij}\}$, and let $T$ play the same role for $\{F^{(2)}_{ij}\}$. We remark in passing that $E$ is generated by $\{E_{ij}^{(3)}\} \cup \{S\} \cup (\CK \oplus \CK)$, and similarly, $F$ is generated by $\{F_{ij}^{(2)}\} \cup \{T\} \cup (\CK \oplus \CK)$. Then
\begin{align*}
\text{ind}_{F^{(2)}_{11}(H_+ \oplus H_-)}&(T) = \text{ind}(\tau_F) = \text{ind}(\tau_E) \\
&= \text{ind}_{E^{(3)}_{11}(H_+ \oplus H_-)}(S) = \text{ind}_{F^{(2)}_{11}(H_+ \oplus H_-)}(WSW^*).
\end{align*}
It follows from \cite[Theorem 3.1]{bdf} that $WSW^*$ and $T$ are unitarily equivalent modulo the compact operators. Thus there is $Z_1 \in U(F^{(2)}_{11}(H_+ \oplus H_-))$ such that $Z_1 WSW^* Z_1^* = T + B$, where $B$ is compact. Let $Z_i = F^{(2)}_{i1} Z_1 F^{(2)}_{1i}$ for $i > 1$, and let $Z = \sum_{i=1}^n Z_i + (1 - \sum_{i=1}^n F^{(2)}_{ii})$. Then
\begin{align*}
Z F^{(2)}_{ij} Z^* &= Z_i F^{(2)}_{ij} Z_j^* = F^{(2)}_{i1} Z_1 F^{(2)}_{1i} F^{(2)}_{ij} F^{(2)}_{j1} Z_1^* F^{(2)}_{1j} \\
&= F^{(2)}_{i1} Z_1 F^{(2)}_{11} Z_1^* F^{(2)}_{1j} = F^{(2)}_{i1} Z_1 Z_1^* F^{(2)}_{1j} = F^{(2)}_{i1} F^{(2)}_{11} F^{(2)}_{1j} = F^{(2)}_{ij}, \\
\noalign{and}
Z WSW^* Z^* &= Z W E_{11}^{(3)} S E^{(3)}_{11} W^* Z^* = Z F^{(2)}_{11} WSW^* F^{(2)}_{11} Z^* \\
&= Z_1 WSW^* Z_1^* = T + B.
\end{align*}
Therefore $Z W E W^* Z^* = F$, and hence $E$ and $F$ are isomorphic.
\end{proof}

We summarize the results of this section:

\begin{Theorem} \label{thm isomorphism invariant}

Let $E$ and $F$ be the algebras in the center places of two versions of the extension \eqref{eqn the extension} (with corresponding Busby invariants $\tau_E$ and $\tau_F$). Suppose that $\emph{ind}(\tau_E)$ and $\emph{ind}(\tau_F)$ have the same orbit under the action of $\eta$ and $\eps$. Then after adjusting with $\eta$ and/or $\eps$ if necessary so that the indices agree, we have that $E$ is isomorphic to $F$ if and only if $\overline{m}(\{E_{ij}\}) = \overline{m}(\{F_{ij}\})$ for any families of matrix units in $E$ and $F$ that lift $\{e_{ij} \otimes 1\}$.

\end{Theorem}

\begin{Definition}

Because of this theorem, we may write $\overline{m}(E)$ instead of $\overline{m}(\{E_{ij}\})$.

\end{Definition}

\section{The invariant $\overline{m}$} \label{sec the invariant}

For convenience, in this section we will view elements of $\IZ^2$ as columns. Thus the index ind$(\tau)$ of the extension \eqref{eqn the extension} will be written $a = (a_+,a_-)^T$. We will describe the invariant $\overline{m}(E)$ from Section \ref{sec extensions}.  It lies in the group $\IZ^2 / (\IZ a + n \IZ^2)$. Let $c = \text{gcd}(a_+,a_-)$, and let $a' \equiv (a_+',a_-')^T = c^{-1} (a_+,a_-)^T$. Choose $b = (b_+,b_-)^T \in \IZ^2$ with $-a_+'  b_- + a_-'  b_+ = 1$. Then $\{a',b\}$ is a basis for $\IZ^2$. Let $d = \text{gcd}(c,n) = \text{gcd}(a_+,a_-,n)$.

\begin{Lemma} \label{lem Za + nZ^2}

$\IZ a + n \IZ^2 = d \IZ a' + n \IZ b$.

\end{Lemma}

\begin{proof}
First note that $\IZ a + n \IZ a' = (c \IZ + n \IZ) a' = d \IZ a'$. Now we have
\begin{equation*}
\IZ a + n \IZ^2 = \IZ a + n \IZ a' + n \IZ b = d \IZ a' + n \IZ b. \qedhere
\end{equation*}
\end{proof}

\begin{Proposition} \label{prop D}

$\IZ^2 / (\IZ a + n \IZ^2) \cong \frac{\IZ}{d \IZ} \oplus \frac{\IZ}{n \IZ}$. The isomorphism is given as follows. Let $A = \begin{psmallmatrix} 0 & -1 \\ 1 & \phantom{-} 0 \end{psmallmatrix}$. For $k = (k_1,k_2)^T \in \IZ^2$,
\begin{equation*}
k + \IZ a + n \IZ^2 \mapsto (k^T A b + d \IZ, k^T A^T a' + n \IZ).
\end{equation*}

\end{Proposition}

\begin{proof}
Note that $\IZ^2 = \IZ a' + \IZ b$ is an internal direct sum. Then Lemma \ref{lem Za + nZ^2} implies that
\begin{equation*}
\IZ^2 / (\IZ a + n \IZ^2) =  \frac{\IZ a' + \IZ b}{d \IZ a' + n \IZ b} = \frac{\IZ}{d \IZ} a' + \frac{\IZ}{n \IZ} b \cong \frac{\IZ}{d \IZ} \oplus \frac{\IZ}{n \IZ}.
\end{equation*}
Note that
\begin{align*}
(a')^T A b &= 1 \\
(a')^T A^T a' &= 0 \\
b^T A b &= 0 \\
b^T A^T a' &= 1.
\end{align*}
Then for any $k = \lambda a' + \mu b \in \IZ^2$, we have $(k^T A b + d \IZ, k^T A^T a' + n \IZ) = (\lambda + d \IZ, \mu + n \IZ)$.
\end{proof}

\begin{Remark} \label{rem main case}

We are mainly interested in the case $a = (-1,1)^T$. In this case $c = d = 1$, so $\IZ^2 / (\IZ a + n \IZ^2) \cong \frac{\IZ}{n \IZ}$. For $k \in \IZ^2$,
\begin{equation*}
k + \IZ a + n \IZ^2 \mapsto k^T A^T a = k_+ + k_-.
\end{equation*}

\end{Remark}

\begin{Corollary} \label{rem index (-1,1)}

Let $E$ be the center algebra of an instance of the extension \eqref{eqn the extension} with index $(-1,1)$. Let $\{E_{ij}\}$ be matrix units in $E$ lifting $\{e_{ij} \otimes 1\}$, and let $k = (k_+,k_-)$ be the defects of $\{E_{ij}\}$. Then $\overline{m}(E) \equiv k_+ + k_- \pmod n$.

\end{Corollary}

Theorem \ref{thm isomorphism invariant} may be adapted to this case as follows.

\begin{Theorem} \label{thm invariant for index (-1,1)}

Let
\begin{align*}
0 \to \CK \oplus \CK \to &E \to M_n \otimes C(\IT) \to 0 \\
0 \to \CK \oplus \CK \to &F \to M_p \otimes C(\IT) \to 0
\end{align*}
be two instances of \eqref{eqn the extension} with index $(-1,1)$. Then $E$ and $F$ are isomorphic if and only if $n=p$ and $\overline{m}(E) = \overline{m}(F)$.

\end{Theorem}

\begin{Proposition} \label{prop tensor with matrices}

Let the extension \eqref{eqn the extension} be given, with index $(-1,1)$. Suppose that $t \in \IN$ is a common divisor of $\overline{m}(E)$ and $n$ (where we let $\overline{m}(E) \in \{ 0,1,2,\ldots,n-1 \}$). Let $\overline{m}(E) = t \ell$ and $n = t p$, and let
\begin{equation*}
0 \to \CK \oplus \CK \to F \to M_p \otimes C(\IT) \to 0
\end{equation*}
be another extension having index $(-1,1)$, and with $\overline{m}(F) = \ell$. ($F$ is uniquely determined, by Theorem \ref{thm isomorphism invariant}.) Then $E \cong M_t \otimes F$.

Conversely, if $E$ is isomorphic to $M_t \otimes F$ for some $C^*$-algebra $F$, then $t$ divides $n$, and letting $n = tp$ there is an extension
\begin{equation*}
0 \to \CK \otimes \CK \to F \to M_p \otimes C(\IT) \to 0
\end{equation*}
with index $(-1,1)$ and $\overline{m}(E) = t \overline{m}(F)$.

\end{Proposition}

\begin{proof}
Using Lemma \ref{lem shift defect by index} we may assume that $k = (t\ell,0)$. Then, as in the proof of that lemma, we may write
\begin{align*}
H_+ &= \bigl[ \IC^{tp} \otimes ( \ell^2 \oplus L_+ ) \bigr] \oplus \IC^{t\ell} \\
&\cong \IC^t \otimes \Bigl( \bigl[ \IC^p \otimes (\ell^2 \oplus L_+) \bigr] \oplus \IC^\ell \Bigr) \\
H_- &= \IC^{tp} \otimes ( \ell^2 \oplus L_- ) \\
&\cong \IC^t \otimes \bigl( \IC^p \otimes (\ell^2 \oplus L_-) \bigr) \\
 S_+ &= e_{11} \otimes \Bigl( \bigl[ e_{11} \otimes (R \oplus V_+) \bigr] \oplus 0 \Bigr) \\
 S_- &= e_{11} \otimes \bigl( e_{11} \otimes (R^* \oplus V_-) \bigr) \\
E_{(i-1)t + k, (j-1)t  + \ell}^+ &= e_{ij} \otimes \Bigl( \bigl[ e_{k \ell} \otimes (1 \oplus 1) \bigr] \oplus 0 \Bigr) \\
E_{(i-1)t + k, (j-1)t  + \ell}^- &= e_{ij} \otimes \bigl( e_{k \ell} \otimes (1 \oplus 1) \bigr),
\end{align*}
where in the last two lines, $1 \le i,j \le t$ and $1 \le k,\ell \le p$. The result is now clear.

For the converse, suppose that $E = M_t \otimes F$. It is easy to see that $F$ must be unital. Let $\{G_{ij} : 1 \le i,j \le t\} \subseteq E$ be matrix units for $M_t \otimes 1_F$, so we may identify $F$ with $G_{11} E G_{11}$. Then $\sum_{i=1}^t \pi(G_{ii}) = \sum_{j=1}^n e_{jj} \otimes 1$. Since $\{\pi(G_{ii})\}$ are pairwise equivalent it must be the case that $t$ divides $n$, and that $\pi(G_{11})$ is Murray-von Neumann equivalent to $\sum_{j=1}^{t^{-1}n} e_{jj} \otimes 1$. Let $n = tp$. Then $\pi(F) \cong M_p \otimes C(\IT)$, and $\ker(\pi|_F) = \CK(G_{11}^+ H_+) \oplus \CK(G_{11}^- H_0)$. Thus we have  the extension
\begin{equation*}
0 \to \CK \oplus \CK \to F \to M_p \otimes C(\IT) \to 0.
\end{equation*}
The index is calculated from the same unitary element $e_{11} \otimes z$, hence is still $(-1,1)$. Since $\sum_{i=1}^t G_{ii} = 1$ and the $\{G_{ii}\}$ are pairwise equivalent, the defect of $E$ equals $t$ times the defect of $F$. Therefore $t$ divides $\overline{m}(G)$.
\end{proof}

\section{The $C^*$-algebras for continued fractions of rational numbers}
\label{sec continued fractions of rationals}

In \cite{ms} the Effros-Shen algebras are presented as $C^*$-algebras of ample \'etale groupoids constructed from certain categories of paths. Categories of paths were introduced in \cite{spi} as a generalization of directed graphs, higher rank graphs, and submonoids of groups, and it was shown that the construction of $C^*$-algebras for these special examples works generally in the setting of categories of paths. The Effros-Shen algebras are AF $C^*$-algebras associated to irrational numbers by means of their simple continued fraction expansions. The construction in \cite{ms} proceeds from an alternative nonsimple continued fraction expansion, and the groupoids that result have quite different properties from the AF groupoids of the usual construction. In this section we consider a variation on the examples of \cite{ms} that correspond to rational numbers. The construction is closely related to the one in \cite{ms}, and in fact occurs there as a building block in that paper. We will use the details of the construction only in a few places (Definition \ref{def Phi and Psi} and Lemmas \ref{lem phi and psi} and \ref{lem dimension and defect of examples}). We will now summarize the basic information needed, and we refer to \cite[sections 3,4]{ms} for details. In particular, we will not give background on categories of paths and the construction of groupoids and $C^*$-algebras from them.

The construction in \cite{ms} begins with a sequence $k = (k_i)_{i=1}^\infty \in \prod_1^\infty \IN$ with the property that $k_i \not= 0$ for infinitely many $i$. The category of paths $\Lambda \equiv \Lambda(k)$ can be thought of as a directed graph with \textit{identifications} (or \textit{relations}); we use the usual notation for directed graphs (or small categories), e.g. $r$ and $s$ for range and source, etc. Here is the sketch:
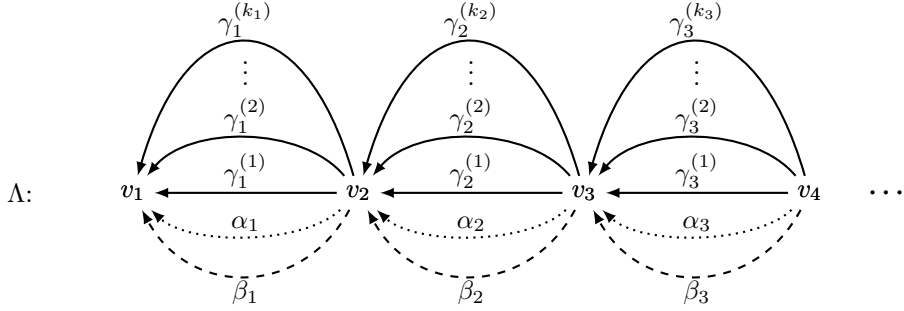
\[
\begin{tikzpicture}[xscale=3,yscale=2]

\node (00) at (.5,0) [rectangle] {$\Lambda$:};

\node (10) at (1,0) [rectangle] {$v_1$};
\node (20) at (2,0) [rectangle] {$v_2$};
\node (30) at (3,0) [rectangle] {$v_3$};
\node (40) at (4,0) [rectangle] {$v_4$};

\node (4half0) at (4.35,0) [rectangle] {$\boldsymbol{\cdots}$};
\node (1half1) at (1.5,.85) [rectangle] {$\boldsymbol{\vdots}$};
\node (2half1) at (2.5,.85) [rectangle] {$\boldsymbol{\vdots}$};
\node (23half1) at (3.5,.85) [rectangle] {$\boldsymbol{\vdots}$};

\draw[-latex,thick] (20) -- (10) node[pos=0.5,inner sep=0.5pt,above=1pt] {$\gamma^{(1)}_1$};
\draw[-latex,thick] (20) .. controls (1.75,.45) and (1.25,.45) .. (10) node[pos=0.5,inner sep=0.5pt,above=1pt] {$\gamma^{(2)}_1$};
\draw[-latex,thick] (20) .. controls (1.75,1.3) and (1.25,1.3) .. (10) node[pos=0.5,inner sep=0.5pt,above=1pt] {$\gamma^{(k_1)}_1$};

\draw[-latex,thick] (30) -- (20) node[pos=0.5,inner sep=0.5pt,above=1pt] {$\gamma^{(1)}_2$};
\draw[-latex,thick] (30) .. controls (2.75,.45) and (2.25,.45) .. (20) node[pos=0.5,inner sep=0.5pt,above=1pt] {$\gamma^{(2)}_2$};
\draw[-latex,thick] (30) .. controls (2.75,1.3) and (2.25,1.3) .. (20) node[pos=0.5,inner sep=0.5pt,above=1pt] {$\gamma^{(k_2)}_2$};

\draw[-latex,thick] (40) -- (30) node[pos=0.5,inner sep=0.5pt,above=1pt] {$\gamma^{(1)}_3$};
\draw[-latex,thick] (40) .. controls (3.75,.45) and (3.25,.45) .. (30) node[pos=0.5,inner sep=0.5pt,above=1pt] {$\gamma^{(2)}_3$};
\draw[-latex,thick] (40) .. controls (3.75,1.3) and (3.25,1.3) .. (30) node[pos=0.5,inner sep=0.5pt,above=1pt] {$\gamma^{(k_3)}_3$};


\node (10) at (1,0) [rectangle] {$v_1$};
\node (20) at (2,0) [rectangle] {$v_2$};
\node (30) at (3,0) [rectangle] {$v_3$};
\node (40) at (4,0) [rectangle] {$v_4$};

\draw[-latex,thick,dotted] (20) .. controls (1.75,-.35) and (1.25,-.35) .. (10) node[pos=0.5,inner sep=0.5pt,above=1pt] {$\alpha_1$};
\draw[-latex,thick,dotted] (30) .. controls (2.75,-.35) and (2.25,-.35) .. (20) node[pos=0.5,inner sep=0.5pt,above=1pt] {$\alpha_2$};
\draw[-latex,thick,dotted] (40) .. controls (3.75,-.35) and (3.25,-.35) .. (30) node[pos=0.5,inner sep=0.5pt,above=1pt] {$\alpha_3$};

\draw[-latex,thick,dashed] (20) .. controls (1.75,-.7) and (1.25,-.7) .. (10) node[pos=0.5,inner sep=0.5pt,below=1pt] {$\beta_1$};
\draw[-latex,thick,dashed] (30) .. controls (2.75,-.7) and (2.25,-.7) .. (20) node[pos=0.5,inner sep=0.5pt,below=1pt] {$\beta_2$};
\draw[-latex,thick,dashed] (40) .. controls (3.75,-.7) and (3.25,-.7) .. (30) node[pos=0.5,inner sep=0.5pt,below=1pt] {$\beta_3$};

\end{tikzpicture}
\]
Thus $s(\alpha_i) = s(\beta_i) = s(\gamma_i^j) = v_{i+1}$ and $r(\alpha_i) = r(\beta_i) = r(\gamma_i^j) = v_i$ for all $i$ and all $1 \le j \le k_i$. The identifications are: $\alpha_i \beta_{i+1} = \beta_i \alpha_{i+1}$, $i \ge 1$. Thus we may think of the edges $\alpha,\beta$ as ``commuting'' with each other, while the edges $\gamma$ can be thought of as ``walls'' that block interaction between $\alpha$ and $\beta$ edges on their sides. Thus an element of $\Lambda$ that does not involve any $\gamma$ edges can be written uniquely in the form $\alpha_h \alpha_{h+1} \cdots \alpha_{h + i - 1} \beta_{h + i} \beta_{h + i + 1} \cdots \beta_{h + i + j - 1}$, which we abbreviate as $v_h \alpha^i \beta^j$, or just $\alpha^i \beta^j$ if the range is understood. We let $\Lambda_1$ be the subcategory generated by the edges $\{ \alpha_i,\beta_i : i \ge 1 \}$, and $\Lambda_2$ the subcategory generated by the edges $\{ \gamma_i^{j} : i \ge 1, 1 \le j \le k_i \}$. Thus $\Lambda_1$ is a 2-graph, and $\Lambda_2$ is a directed graph. ($\Lambda$ is an \textit{amalgamation} of $\Lambda_1$ and $\Lambda_2$ as in \cite[section 11]{spi}.) In order to build the groupoid $G$ we must first describe the unit space $X \equiv G^{(0)} = v_1 \Lambda^\infty$, which is a compact Hausdorff space. The elements of $X$ are (equivalence classes of) infinite paths in $\Lambda$ with range $v_1$. (The equivalence relation corresponds to the commutativity of $\alpha$ and $\beta$ edges.) Thus we may describe an element of $X$ in one of two forms:
\begin{align*}
\quad& v_1 \alpha^{i_1} \beta^{j_1} \gamma^{(r_1)} \alpha^{i_2} \beta^{j_2} \gamma^{(r_2)} \cdots, \\
\quad& v_1 \alpha^{i_1} \beta^{j_1} \gamma^{(r_1)} \cdots \alpha^{i_m} \beta^{j_m} \gamma^{(r_m)} \alpha^p \beta^q,
\end{align*}
where $i_h, j_h \ge 0$, $r_h \ge 1$, $p,q \ge 0$, and $p + q = \infty$, and where we have omitted subscripts from the edges. (We see from the second of these that the description of elements of $X$ as ``infinite paths'' is not quite accurate. The correct description as \textit{directed hereditary sets} can be found in \cite{spi}.) To define the topology on $X$ we use the \textit{cylinder sets} $Z(\mu) = \{ \mu x : \mu \in v_1 \Lambda, x \in r(\mu) \Lambda^\infty \}$, i.e. all (equivalence classes of) infinite paths that begin with $\mu$. It follows from \cite[Proposition 5.6]{ms} that the sets $\{ Z(\mu) : \mu \in v_1 \Lambda \} \cup \{ Z(\mu) \setminus Z(\mu\lambda) : \mu \in v_1 \Lambda, \lambda \in \Lambda^1 \}$ form a base of compact open sets for a totally disconnected compact Hausdorff topology on $X$. (We write $\Lambda^1$ for the set of edges in $\Lambda$.) The groupoid $G$ consists of equivalence classes $G = \{ [\mu,\nu,x] : \mu,\nu \in v_1 \Lambda, s(\mu) = s(\nu), x \in s(\mu) \Lambda^\infty \}$. The equivalence relation is generated by $(\mu,\nu,\eps y) \sim (\mu\eps,\nu\eps,y)$. Source and range are given by $s([\mu,\nu,x]) = \nu x$, $r([\mu,\nu,x]) = \mu x$, and the product is given by $[\mu,\nu,x] \cdot [\nu,\eta,x] = [\mu,\eta,x]$, where the equivalence relation allows all composable pairs to be written this way. The sets $[\mu,\nu,E] = \{ [\mu,\nu,x] : x \in E \}$, where $E \subseteq X$ is a compact open subset, form a base of compact open bisections for an ample Hausdorff topology on $G$. The nonsimple continued fraction $[0,1,k_1,1,k_2,\ldots]$ converges to an irrational number $\theta \in (0,1)$, and it is shown in \cite[Theorem 7.2]{ms} that  $C^*(G)$ is isomorphic to the Effros-Shen algebra $A_\theta$. The proof relies on a decomposition of $C^*(G)$ as a direct limit of $C^*$-algebras of subgroupoids $G_i$ of $G$. From \cite[Theorem 3.18]{ms} we have that $G_i = \{ [\mu\alpha^a\beta^b,\nu\alpha^c\beta^d,x] \in G : \mu,\nu \in v_1 \Lambda v_{i+1}, 0 \le a,b,c,d < \infty \}$. In words, the first two coordinates do not use $\gamma$ edges after the $i$th spot.

In this section we will consider sequences $k \in \prod_1^\infty \IN$ such that $k_i \not= 0$ for only finitely many $i$. The constructions of $\Lambda$ and $G$ proceed exactly as in \cite{ms}, but of course the results are different. Fix $k \in \prod_1^\infty \IN$ with $k_i \not= 0$ for only finitely many $i$, and let $\Lambda$, $G$ be the corresponding category of paths and Hausdorff ample groupoid. Let $h = \max \{i : k_i \not= 0\}$. In this case $X = F_h$, where $F_h$ is in general a proper closed subset of $X$, defined in \cite[Definition 4.1]{ms}. (For example, if $k' \in \prod_1^\infty \IN$ is any sequence that agrees with $k$ for $i \le h$, but is nonzero infinitely often, then $X(k)$ is the same as $F_h(k')$ as defined in \cite{ms}. Thus $C^*(G) = C^*(G(k')|_{F_h})$, the restriction of the groupoid to the closed subset $F_h$ of its unit space.) It is shown after \cite[Corollary 4.16]{ms} that the algebra $C^*(G)$ is the central algebra of an extension \eqref{eqn the extension} with index $(-1,1)$. We will associate this algebra to the rational number $[0,1,k_1,1,k_2,\ldots,1,k_h]$. In order to describe the size of the matrix algebra in the quotient, we require some definitions (see \cite[Definition 4.9]{ms}).

\begin{Definition} \label{def Phi and Psi}

Let $0 \le f \le h$. We define
%
\begin{align*}
\Psi_f &= \{\mu \in v_1 \Lambda v_{f+1} : \mu_f \in \Lambda_2^1 \} \\
\psi_f &= |\Psi_f| \\
\Phi_f &= \{\mu \in v_1 \Lambda : |\mu| \le f \text{ and } \mu_{|\mu|} \in \Lambda_2^1 \} = \bigsqcup_{\ell=0}^f \Psi_\ell \\
\phi_f &= |\Phi_f| = \sum_{\ell=0}^f \psi_f.
\end{align*}

\end{Definition}

\begin{Lemma} \label{lem phi and psi}
\begin{enumerate}

\item \label{lem phi and psi one} $\psi_f = k_f \sum_{\ell=0}^{f-1} (f-\ell)\psi_\ell$ for $f \ge 1$, and $\psi_0 = 1$.

\item \label{lem phi and psi two} $\phi_f = k_f \sum_{\ell=0}^{f-1} \phi_\ell + \phi_{f-1}$ for $f \ge 1$, and $\phi_0 = 1$.

\item \label{lem phi and psi three} $\phi_f$ and $\sum_{\ell=0}^{f-1} \phi_\ell$ are relatively prime.

\end{enumerate}
\end{Lemma}

\begin{proof}
It is clear that $\Psi_0 = \Phi_0 = \{v_1\}$, and hence $\psi_0 = \phi_0 = 1$. Let $\mu \in \Psi_f$. Then $\mu = \mu' \gamma_f^{(r)}$ for some $1 \le r \le k_f$. Then there is $\ell < f$, $\nu \in \Psi_\ell$, and $i,j \ge 0$ with $i + j = f - \ell - 1$ such that $\mu' = \nu \alpha^i \beta^j$. There are $\psi_\ell$ choices for $\nu$, and $f - \ell$ choices for $(i,j)$. Thus $\psi_f = k_f \sum_{\ell=0}^{f-1} (f - \ell) \psi_\ell$. This proves \eqref{lem phi and psi one}. To prove  \eqref{lem phi and psi two}, first note that
\begin{equation} \label{eqn phi and psi}
\sum_{\ell=0}^{f-1} (f-\ell) \psi_\ell
= \sum_{\ell=0}^{f-1} \sum_{p=0}^{f - \ell - 1} \psi_\ell
= \sum_{p=0}^{f-1} \sum_{\ell=0}^{f-p-1} \psi_\ell
= \sum_{p=0}^{f-1} \phi_{f-p-1}.
\end{equation}
Now we have
\begin{align*}
\phi_f &= \sum_{\ell=0}^f \psi_\ell
= \psi_f + \sum_{\ell=0}^{f-1} \psi_\ell
= k_f \sum_{\ell=0}^{f-1}(f - \ell) \psi_\ell + \phi_{f-1}, \text{ by \eqref{lem phi and psi one},} \\
&= k_f \sum_{\ell=0}^{f-1} \phi_{f-\ell-1} + \phi_{f-1}
= k_f \sum_{\ell=0}^{f-1} \phi_\ell + \phi_{f-1}, \text{ by \eqref{eqn phi and psi}.} 
\end{align*}
To prove \eqref{lem phi and psi three}, let $c$ be the greatest common divisor of $\phi_f$ and $\sum_{\ell=0}^{f-1} \phi_\ell$. Then \eqref{lem phi and psi two} implies that $c$ divides $\phi_{f-1}$. It follows that $c$ divides $\sum_{\ell_0}^{f-2} \phi_\ell$. Now \eqref{lem phi and psi two}, applied to $f-1$, implies that $c$ divides $\phi_{f-2}$. Continuing this process, we see that $c$ divides $\phi_0 = 1$.
\end{proof}

\begin{Remark}

It follows from Lemma \ref{lem phi and psi}\eqref{lem phi and psi two} that $\sum_{\ell=0}^{h-1} \phi_\ell < \phi_h$.

\end{Remark}

\begin{Lemma} \label{lem dimension and defect of examples}

The size of the matrices in the quotient $C^*(G)/(\CK \oplus \CK)$ is $\phi_h$. The sum of the defects is $\overline{m}(C^*(G)) = \sum_{\ell=0}^{h-1} \phi_\ell$ (as in Remark \ref{rem main case}).

\end{Lemma}

\begin{proof}
The first statement is given by \cite[Corollary 4.12]{ms}. For the second statement we analyze more closely the extension \eqref{eqn the extension}
\begin{equation*}
0 \to \CK(H_+) \oplus \CK(H_-) \to C^*(G) \xrightarrow{\pi} M_{\phi_h} \otimes C(\IT) \to 0,
\end{equation*}
where $M_{\phi_h} \otimes C(\IT) = C^*(G|_{F_h^\infty})$ in the notation of \cite[section 4]{ms}. (Here, $F_h^\infty = \{\mu \alpha^\infty \beta^\infty : \mu \in \Phi_h\}$.) By \cite[Corollary 4.12]{ms}, $e_{\mu,v_1} \otimes 1 = \chi_{[\mu,\alpha^{|\mu|}, \alpha^\infty \beta^\infty]}$ for $\mu \in \Phi_h$. (Here, we are letting $e_{v_1,v_1}$ play the role of $e_{11} \in M_n$; $v_1$ is the unique element of $\Phi_h$ with length zero.) We lift these to matrix units in $C^*(G)$ by
\begin{align*}
E_{\mu,v_1} &= \chi_{[\mu \alpha^{h - |\mu|}, \alpha^h, Z(v_{h+1})]}, \text{ for } \mu \in \Phi_h \setminus \{v_1\}. \\
\noalign{Then}
E_{\mu,v_1}^* E_{\mu,v_1} &= \chi_{Z(\alpha^h)} \\
E_{\mu,v_1} E_{\mu,v_1}^* &= \chi_{Z(\mu \alpha^{h - |\mu|})} \\
E_{\mu,v_1}^2 &= \chi_{ [ \mu \alpha^{h - |\mu|}, \alpha^h, Z(v_{h+1}) ] \cdot [\mu \alpha^{h - |\mu|}, \alpha^h, Z(v_{h+1})] } \\
&= 0,
\end{align*}
since $Z(\alpha^h) \cap Z(\mu \alpha^{h - |\mu|}) = \varnothing$. Moreover, for $\mu \not= \nu \in \Phi_h$, $E_{\mu,v_1}^* E_{\nu,v_1} = 0$, since $Z(\mu \alpha^{h - |\mu|}) \cap Z(\nu \alpha^{h - |\nu|}) = \varnothing$. Therefore $\{ E_{\mu,v_1} : \mu \in \Phi_h \setminus \{v_1\} \}$ define a family of $\phi_h \times \phi_h$ matrix units in $C^*(G)$, and $\pi(E_{\mu,v_1}) = e_{\mu,v_1} \otimes 1$.

We now calculate $m(\{E_{\mu,\nu} : \mu,\nu \in \Phi_h\})$.  $H_+$ has orthonormal basis $\{ \delta_{\mu \alpha^\infty \beta^i} : \mu \in \Phi_h, i \in \IN \}$, and $H_-$ has orthonormal basis $\{ \delta_{\mu \alpha^i \beta^\infty} : \mu \in \Phi_h, i \in \IN \}$.  Since $E_{\mu \mu} = \chi_{Z(\mu \alpha^{h - |\mu|})}$, we have
\begin{align*}
E_{\mu \mu}^+ H_+ &= \text{span} \{ \delta_{ \mu \alpha^\infty \beta^i} : i \ge 0 \} \\
E_{\mu \mu}^- H_- &= \text{span} \{\delta_{\mu \alpha^i \beta^\infty} : i \ge h - |\mu| \}. \\
\noalign{Then}
\sum_{\mu \in \Phi_h} E_{\mu \mu}^+ &= 1_{H_+} \\
(1_{H_-} - \sum_{\mu \in \Phi_h} E_{\mu \mu}^- )H_- &= \text{span}\{ \delta_{\mu \alpha^i \beta^\infty} : \mu \in \Phi_h, 0 \le i < h - |\mu| \}. \\
\noalign{Therefore}
\text{rank}(1_{H_-} - \sum_{\mu \in \Phi_h} E_{\mu \mu}^-)
&= \sum_{\mu \in \Phi_h} (h - |\mu|) 
= \sum_{\ell=0}^h (h - \ell) \psi_\ell 
= \sum_{\ell=0}^{h-1} (h - \ell) \psi_\ell \\
&= \sum_{p = 0}^{h-1} \phi_{h - p - 1}, \text{ by \eqref{eqn phi and psi},}  \\
&= \sum_{\ell=0}^{h-1} \phi_\ell.
\end{align*}
Thus $m(\{E_{\mu,\nu} : \mu,\nu \in \Phi_h \}) = (0,\sum_{\ell=0}^{h-1} \phi_\ell)$.
\end{proof}

\begin{Lemma} \label{lem E determines k_i}

The assignment $(k_i)_{i=1}^\infty \mapsto C^*(G)$ described at the beginning of this section is a one-to-one correspondence between finitely nonzero sequences in $\IN$ and isomorphism classes of algebras $E$ in the central spot of extensions \eqref{eqn the extension}, having index $(-1,1)$ and such that $\overline{m}(E)$ is relatively prime to $n$, the order of the matrices in the quotient spot of \eqref{eqn the extension} (and where we interpret $\overline{m}(E)$ as an integer, as in Corollary \ref{rem index (-1,1)}). (The exceptional case $\overline{m}(E) = 0$ and $n=1$ corresponds to the sequence identically equal to zero. In this case, the algebra $E$ is isomorphic to $C^*(R \oplus R^*)$, where $R \in B(\ell^2)$ is the right shift operator.)

\end{Lemma}

\begin{proof}
Let $(k_i)_{i\ge1}$ be a finitely nonzero sequence in $\IN$, and let $G$ be the corresponding groupoid. Let $h$ be the largest index for which $k_h > 0$ (and we assume that $(k_i)$ is not the zero sequence). Then $C^*(G)$ is in the center of an extension \eqref{eqn the extension}:
\begin{equation*}
0 \to \CK \oplus \CK \to C^*(G) \to M_{\phi_h} \otimes C(\IT) \to 0,
\end{equation*}
and $\overline{m}(C^*(G)) = \sum_{\ell=0}^{h-1} \phi_\ell$. Since $\phi_h = k_h \sum_{\ell=0}^{h-1} \phi_\ell + \phi_{h-1}$, the Euclidean algorithm allows us to calculate $k_h$ and $\phi_{h-1}$ from $C^*(G)$. Then we may also calculate $\sum_{\ell=0}^{h-2} \phi_\ell = \sum_{\ell=0}^{h-1} \phi_\ell - \phi_{h-1}$. Since $\phi_{h-1} = k_{h-1} \sum_{\ell=0}^{h-2} \phi_\ell + \phi_{h-2}$, another application of the Euclidean algorithm allows us to calculate $k_{h-1}$ and $\phi_{h-2}$. Continuing this process, we see that we may calculate the sequence $(k_i)_{i=0}^h$, hence by appending zeros we obtain $(k_i)_{i=0}^\infty$. Thus the map from sequences in $\IN$ to isomorphism classes of extensions described in the statement is one-to-one.

Now we show that this map is onto. The proof uses a modification of usual calculation of the greatest common divisor of two numbers. Let $0 < m < n$ with $m$ and $n$ relatively prime. By the Euclidean algorithm we may write
\begin{equation*}
n = q_0 m + r_1, \ 0 \le r_1 < m.
\end{equation*}
In fact, $r_1 > 0$ since $m$ does not divide $n$. Then gcd$(r_1,m)$ divides $n$, hence $r_1$ and $m$ are relatively prime. We again use the Euclidean algorithm to write
\begin{equation*}
r_1 = q_1 (m - r_1) + r_2, \ 0 \le r_2 < m - r_1.
\end{equation*}
Note that $r_2 \le r_1$. If $r_2 > 0$ then gcd$(r_2, m - r_1)$ divides $r_1$, hence divides $m$. Therefore $r_2$ and $m - r_1$ are relatively prime. We continue this process to obtain $q_\ell$ and $r_\ell$ with
\begin{equation*}
r_\ell = q_\ell(m - \sum_{i=1}^\ell r_i) + r_{\ell + 1}, \ 0 \le r_{\ell+1} < m - \sum_{i=1}^\ell r_i.
\end{equation*}
While $r_\ell > 0$ we have
\begin{align*}
&m > m - r_1 > m - r_1 - r_2 > \cdots \\
&\text{gcd}(r_\ell,m - \sum_{i=1}^{\ell-1} r_i) = 1.
\end{align*}
Let $h$ be such that $r_h = 0$ and $r_{h-1} > 0$. Then

\begin{align*}
r_{h-1} &= q_{h-1} (m - \sum_{\ell=1}^{h-1} r_\ell) 
= q_{h-1} (m - \sum_{\ell=1}^{h-2} r_\ell) - q_{h-1} r_{h-1}; \\
(q_{h-1} + 1) r_{h-1} &= q_{h-1} (m - \sum_{\ell=1}^{h-2} r_\ell).
\end{align*}
Since $r_{h-1}$ and $m - \sum_{\ell=1}^{h-2} r_\ell$ are relatively prime, it follows that $r_{h-1}$ divides $q_{h-1}$. Therefore $r_{h-1} = q_{h-1}$, and hence $m - \sum_{\ell=1}^{h-1} r_\ell = 1$.

Now define $k_\ell = q_{h-\ell}$ for $2 \le \ell \le h$, and $k_1 = q_{h-1} - 1$ (and of course, $k_\ell = 0$ for $\ell > h$). We will identify $r_\ell$, $m$ and $n$ in terms of the numbers $\phi_\ell$ corresponding to $(k_i)_{i=1}^\infty$.
\begin{align*}
\phi_0 &= 1; \\
r_{h-1} &= q_{h-1} = k_1 + 1 = \phi_1, \text{by Lemma \ref{lem phi and psi}\eqref{lem phi and psi two}}; \\
m - \sum_{\ell=1}^{h-2} r_\ell &= 1 + r_{h-1} = \phi_0 + \phi_1; \\
r_{h-2} &= q_{h-2} (m - \sum_{\ell=1}^{h-2} r_\ell) + r_{h-1} 
= k_2(\phi_0 + \phi_1) + \phi_1 = \phi_2, \text{ by Lemma \ref{lem phi and psi two}\eqref{lem phi and psi two}}; \\
m - \sum_{\ell=1}^{h-3} r_\ell &= m - \sum_{\ell=1}^{h-2} r_\ell + r_{h-2} 
= \phi_0 + \phi_1 + \phi_2. \\
\noalign{Inductively, we have}
m &= (m - r_1) + r_1 
= \sum_{i=0}^{h-2} \phi_i + \phi_{h-1} = \sum_{i=0}^{h-1} \phi_i; \\
r_1 &= q_1(m - r_1) + r_2 
= k_{h-1} \sum_{i=0}^{h-2} \phi_i + \phi_{h-2} 
= \phi_{h-1}, \text{ by Lemma \ref{lem phi and psi two}\eqref{lem phi and psi two}}; \\
n &= q_0 m + r_1 
= k_h \sum_{i=0}^{h-1} \phi_i + \phi_{h-1} = \phi_h, \text{ by Lemma \ref{lem phi and psi two}\eqref{lem phi and psi two}}. \qedhere
\end{align*}
\end{proof}

\section{Infinite and finite continued fractions}
\label{sec continued fractions}

A \textit{continued fraction} is determined by a sequence $(a_i)_{i=0}^\infty$ with $a_i \in \IZ$ and $a_i \ge 0$ for $i > 0$. The continued fraction is called \textit{simple} if $a_i > 0$ for $i > 0$. An extended real number (i.e. in $\IR \cup \{\infty\}$) can (often) be associated to a continued fraction as follows. For each $n$ we let
\begin{equation*}
[a_0,a_1,\ldots,a_n] = a_0 + \cfrac{1}{a_1 +
\cfrac{1}{a_2 + \dotsb
\cfrac{1}{a_n}}}
\end{equation*}
If some of the $a_n$ equal 0, we use the usual partially defined arithmetic in $\widetilde{\IR} = \IR \cup \{\infty\} \subseteq \IC \cup \{\infty\} = \widetilde{\IC}$, as described after \cite[Remarks 6.2]{ms}. The infinite continued fraction is defined by $[a_0,a_1,\ldots] = \lim_{n\to\infty} [a_0,a_1,\ldots,a_n]$, if the limit exists. If the continued fraction is simple then the limit does exist. It is classical that each irrational number can be represented as an infinite simple continued fraction in exactly one way, while each rational number can be represented as a finite simple continued fraction in exactly two ways ($[a_0,a_1,\ldots,a_{n-1}, a_n + 1] = [a_0,a_1,\ldots,a_n,1]$); see \cite[Chapter 10]{hw}.

\begin{Notation} \label{not repeated string}

For a string $(a_1,\ldots,a_k)$ and $n \in \IN \setminus \{0\}$ we will write $(a_1,\ldots,a_k)^n$ for the string of length $nk$ obtained by concatenating $(a_1,\ldots,a_k)$ with itself $n$ times.

\end{Notation}

For nonsimple continued fractions we have the following.

\begin{Lemma} \label{lem nonsimple continued fraction}

 (\cite[Lemma 3.1]{es})
Let $(a_n)_{n=0}^\infty$ be a sequence  with $a_n \in \IZ$ and  $a_n \ge 0$ for $n > 0$. Suppose that $a_1 > 0$, that $a_{2n} \not= 0$ for infinitely many $n$, and that $a_{2n+1} \not= 0$ for infinitely many $n$. Then the infinite continued fraction $[a_0,a_1,\ldots]$ converges (to an irrational number). 

\end{Lemma}

We will need the following additional case of nonsimple continued fractions.

\begin{Lemma} \label{lem rational continued fraction}

Let $(a_n)_{n=0}^\infty$ be a sequence with $a_n \in \IZ$ and $a_n \ge 0$ for $n > 0$. Let $m \in \IN$, and suppose that $a_{m + 2i - 1} = 1$, $a_{m + 2i} = 0$ for $i \ge 1$. Then $\lim_{n\to\infty} [a_0,a_1,\ldots,a_n] = [a_0,a_1,\ldots,a_m]$ in $\widetilde{\IR}$.

\end{Lemma}

\begin{proof}
We will manipulate finite continued fractions using the partially defined arithmetic of $\widetilde{\IR}$. We have, for $x \in \widetilde{\IR}$,
\begin{align*}
[1,0] &= 1 + \frac{1}{0} = 1 + \infty = \infty \\
[0,0,x] &= 0 + \cfrac{1}{0 + \cfrac{1}{x}} = x \\
[1,0,x] &= 1 + [0,0,x] = 1 + x \\
[(1,0)^n] &= [1,0,(1,0)^{n-1}] = 1 + [0,0,(1,0)^{n-1}] \\
&= 1 + [(1,0)^{n-1}] = 2 + [(1,0)^{n-2}] \\
&= \cdots = n - 1 + [1,0] = n - 1 + \infty = \infty \\
[1,(0,1)^{n-1}] &= 1 + [0,(0,1)^{n-1}] = 1 + [0,0,1,(0,1)^{n-2}] \\
&= 1 + [1,(0,1)^{n-2}] = 2 + [1,(0,1)^{n-3}] \\
&= \cdots = n - 1 + [1] = n .
\end{align*}
Now we can calculate
\begin{align*}
[x,(1,0)^n] &= [x,[(1,0)^n]] = [x,\infty] = x + \frac{1}{\infty} = x + 0 = x \\
[x,1,(0,1)^{n-1}] &= [x,[1,(0,1)^{n-1}]] = [x,n] = x + \frac{1}{n}.
\end{align*}
Letting $x = [a_0,a_1,\ldots,a_m]$, we obtain the conclusion of the lemma.
\end{proof}

\begin{Proposition} \label{prop arbitrary k_i}

Let $(k_i)_{i=1}^\infty$ be a sequence in $\IN$. The infinite continued fraction $[0,1,k_1,1,k_2,1,\ldots]$ converges to a limit $\theta \in \IR$. If $k_i > 0$ infinitely often then $\theta$ is irrational. If, say, $k_n > 0$ and $k_i = 0$ for $i > n$ then $\theta = [0,1,k_1,1,\ldots,1,k_n] \in \IQ$.

\end{Proposition}

\begin{proof}
The first statement follows from Lemma \ref{lem nonsimple continued fraction}, and the second statement follows from Lemma \ref{lem rational continued fraction}.
\end{proof}

The continued fraction in Proposition \ref{prop arbitrary k_i} is realized as a simple continued fraction in the following lemma.

\begin{Lemma} \label{lem k_i to a_n}

Let $(k_i)_{i=1}^\infty$ be a sequence in $\IN$. Let $1 \le p_1 < p_2 < \cdots$ be defined by $\{i : k_i > 0 \} = \{p_1,p_2,\ldots\}$. Then
\begin{equation} \label{eqn k_i to a_n}
[0,1,k_1,1,k_2,1,\ldots] = [0,p_1,k_{p_1},p_2 - p_1, k_{p_2}, p_3 - p_2, k_{p_3}, \ldots],
\end{equation}
where the righthand side of \eqref{eqn k_i to a_n} ends with $k_{p_n}$ if $(p_j)$ is a finite sequence with largest index $n$.

\end{Lemma}

\begin{proof}
We first claim that for $p \ge 1$ we have that $[(0,1)^p,x] = [0,p,x]$. We prove this by induction. It is true when $p = 1$, so we suppose it is true for some $p \ge 1$. Then
\begin{align*}
[(0,1)^{p+1},x]
&= [(0,1)^p,[0,1,x]] = [0,p,[0,1,x]], \text{ by the inductive hypothesis,} \\
&= [0,p,0,[1,x]]
= \cfrac{1}{p + 
\cfrac{1}{0 +
\cfrac{1}{[1,x]}}}
= \frac{1}{p + [1,x]} = \frac{1}{[p+1,x]} = [0,p+1,x].
\end{align*}
We next claim that $[1,(0,1)^{n-1},x] = [n,x]$. We can prove this using the previous claim:
\begin{align*}
[1,(0,1)^{n-1},x] &= 1 + [0,(0,1)^{n-1},x] = 1 + [(0,1)^{n-1},x]^{-1} = 1 +  [0,n-1,x]^{-1} \\
&= 1 + ([n-1,x]^{-1})^{-1} = 1 + [n-1,x] = [n,x].
\end{align*}
Using these claims we get
\begin{align*}
[0,1,k_1,1,k_2,1,\ldots]
&= [0,1,(0,1)^{p_1 - 1}, k_{p_1}, 1, (0,1)^{p_2 - p_1 - 1}, k_{p_2},1, \ldots] \\
&= [(0,1)^{p_1},[k_{p_1},1,(0,1)^{p_2 - p_1 - 1}, k_{p_2},1,\ldots]] \\
&= [0,p_1,[k_{p_1},1,(0,1)^{p_2 - p_1 - 1}, k_{p_2},1,\ldots]] \\
&= [0,p_1,k_{p_1},[1,(0,1)^{p_2 - p_1 - 1}, k_{p_2},1,\ldots]] \\
&= [0,p_1,k_{p_1},[p_2 - p_1, k_{p_2},1,\ldots]] \\
&= \cdots = [0,p_1,k_{p_1},p_2 - p_1,k_{p_2},p_3 - p_2,\ldots].
\end{align*}
The lemma now follows from Proposition \ref{prop arbitrary k_i}.
\end{proof}

\begin{Theorem} \label{thm sequences k_i and [0,1)}

There is a bijective map $(k_i)_{i=1}^\infty \in \prod_1^\infty \IN \mapsto [0,1,k_1,1,k_2,1,\ldots] \in [0,1)$ such that finitely nonzero sequences correspond to rational numbers while infinitely nonzero sequences correspond to irrational numbers. 

\end{Theorem}

\begin{proof}
Let $\theta \in [0,1) \cap \IQ$. Write $\theta$ as a finite simple continued fraction with last term having an even index: $\theta = [0,a_1,a_2,\ldots,a_{2n}]$. Define $p_j = a_1 + a_3 + \cdots + a_{2j-1}$ and $k_{p_j} = a_{2j}$ for $1 \le j \le n$, and set $k_i = 0$ for $i \not\in \{p_j : 1 \le j \le n\}$. Then by Lemma \ref{lem k_i to a_n} we have
\begin{align*}
[0,1,k_1,1,k_2,1,\ldots] &= [0,p_1,k_{p_1},p_2 - p_1,k_{p_2},\ldots,p_n - p_{n-1}, k_{p_n}] \\
&= [0,a_1,a_2,\ldots,a_{2n}] = \theta.
\end{align*}
For $\theta \in (0,1) \setminus \IQ$ we may write $\theta = [0,a_1,a_2,\ldots]$ as an infinite simple continued fraction. Define $p_j = a_1 + a_3 + \cdots + a_{2j-1}$ and $k_{p_j} = a_{2j}$ for $j \ge 1$, and set $k_i = 0$ for $i \not\in \{p_j : j \ge 1 \}$. Again by Lemma \ref{lem k_i to a_n} we have
\begin{align*}
[0,1,k_1,1,k_2,1,\ldots] &= [0,p_1,k_{p_1},p_2 - p_1,k_{p_2},p_3 - p_2, k_{p_3},\ldots] \\
&= [0,a_1,a_2,\ldots] = \theta. \qedhere
\end{align*}
\end{proof}

\begin{Corollary} \label{cor from a rational to an extension}

There is a bijective correspondence between $[0,1) \cap \IQ$ and the isomorphism classes of extension $C^*$-algebras $E$ in \eqref{eqn the extension} having index $(-1,1)$ such that $\overline{m}(E)$ is relatively prime to $n$.

\end{Corollary}

\begin{proof}
This follows from Theorem \ref{thm sequences k_i and [0,1)} and Lemma \ref{lem E determines k_i}.
\end{proof}

Let us describe this correspondence. If a rational number in $(0,1)$ is given (in lowest terms) as $\tfrac{p}{q}$ then we know that $q$ is the size of the matrix algebra in the final term of \eqref{eqn the extension}. In order to find the defect of the extension we may first write $\tfrac{p}{q}$ as a finite simple continued fraction of even length, then use Theorem \ref{thm sequences k_i and [0,1)} to write it in the form $[0,1,k_1,1,\ldots,1,k_h]$. From this sequence $(k_i)$ we may then use Lemma \ref{lem phi and psi}\eqref{lem phi and psi two} to construct the defect $\sum_{\ell=0}^{h-1} \phi_\ell$. These two numbers determine the extension.

For the reverse direction, it was noted at the beginning of the proof of Lemma \ref{lem E determines k_i} how to use the order $n$ of the matrix algebra and the defect in \eqref{eqn the extension} to calculate a sequence $(k_i)$ that is finitely nonzero. From there one obtains the rational number as $[0,1,k_1,1,\ldots,1,k_h]$ (where $h$ is the largest index of a nonzero term).

Note that in the rational case, only one of the two possible simple continued fractions representing $\theta$ is used in this process. It was shown in \cite[Theorem 7.2]{ms} that there is a bijective map from $(0,1) \setminus \IQ$, in the form given by Theorem \ref{thm sequences k_i and [0,1)}, to the class of all Effros-Shen algebras. Together with the results of this paper, we obtain a bijection between $[0,1)$ and a family of related $C^*$-algebras interpolating between the Effros-Shen algebras. In future work we will show that this family gives a strongly continuous $C^*$-bundle.

Let us contrast the above with the usual construction of the Effros-Shen algebras as AF algebras (\cite{es}). Let $\theta \in (0,1) \setminus \IQ$ have simple continued fraction expansion $[0,a_1,a_2,\ldots]$. Define (\cite[Chapter 10]{hw})
\begin{align*}
p_0 &= 0,\ p_1 = 1,\ p_n = a_n p_{n-1} + p_{n-2} \text{ for } n \ge 2 \\
q_0 &= 1,\ q_1 = a_1,\ q_n = a_n q_{n-1} + q_{n-2} \text{ for } n \ge 2.
\end{align*}
Then $\tfrac{p_n}{q_n} = [0,a_1,\ldots,a_n] \to \theta$ as $n \to \infty$. Let $A_n = M_{q_n} \oplus M_{q_{n-1}}$, $T_n = \begin{psmallmatrix} a_n & 1 \\ 1 & 0 \end{psmallmatrix}$, and let $f_n : A_n \hookrightarrow A_{n+1}$ be a $*$-homomorphism with multiplicities given by $T_n$. Then $A_\theta := \underset{\to}{\lim} (A_n,f_n)$ is the \textit{Effros-Shen} AF algebra. If $\theta \in (0,1) \cap \IQ$, and $[0,a_1,\ldots,a_n]$ is one of the two simple continued fraction representations of $\theta$, then the sequences $(p_i)$ and $(q_i)$ are finite with largest index $n$. It seems reasonable to let $M_{q_n} \oplus M_{q_{n-1}}$ play the role of the Effros-Shen algebra for $\theta$, and we see that there are two such natural candidates.

The way that continued fractions behaves implies that for $x > 0$,
\begin{align*}
[0,a_1,\ldots,a_{2n}] &< [0,a_1,\ldots,a_{2n},x] \\
[0,a_1,\ldots,a_{2m-1},x] &< [0,a_1,\ldots,a_{2m-1}].
\end{align*}
Thus in order to make a $C^*$-bundle by interpolating the finite dimensional analogs of Effros-Shen algebra at the rational numbers, the interval $[0,1)$ should be broken by replacing each rational point $\theta$ with two points, $\theta_\pm = \lim_{\eta \to \theta^\pm} \eta$. The two finite dimensional $C^*$-algebras for $\theta$ will be the fibers at the pair of points $\theta_\pm$. In this way we construct a continuous $C^*$-bundle over a disconnected version of $[0,1)$. This will be made precise in future work.

\end{document}